\documentclass[12pt]{article}
\usepackage{amssymb,latexsym,amsmath}
\usepackage{amsthm}
\usepackage{geometry}
\usepackage{graphicx,color}
\usepackage{lineno}
\newtheorem{thm}{Theorem}

\theoremstyle{plain}

\theoremstyle{definition}

\numberwithin {equation}{section}
\usepackage{longtable}

\usepackage{float}
\restylefloat{table}

\begin{document}
\title
{Modified Potra-Pt{\'a}k method to determine the multiple zeros of
nonlinear equations}
\author{ S. Sharifi$^a$\thanks{s.sharifi@iauh.ac.ir}\and M. Ferrara$^b$\thanks{massimiliano.ferrara@unirc.it} \and
 N.M.A. Nik Long$^c$\thanks{nmasri@upm.edu.my}  \and M. Salimi$^d$\thanks{mehdi.salimi@tu-dresden.de}}
\date{}
\maketitle
\begin{center}
$^{a}$Young Researchers and Elite Club, Hamedan Branch, Islamic
Azad University, Hamedan, Iran\\
$^{b}$Department of Law and Economics, University Mediterranea of
Reggio Calabria, Italy\\
$^{c}$Department of Mathematics and Institute for Mathematical Research, Universiti Putra Malaysia, 43400 Serdang, Selangor, Malaysia\\
$^{d}$Center for Dynamics, Department of Mathematics, Technische
Universit{\"a}t Dresden, Germany\\

\end{center}
\maketitle

\maketitle
\begin{abstract}
\noindent

In this paper, we present a third-order iterative method based on
Potra-Pt{\'a}k method to compute the approximate multiple roots of
nonlinear equations. The method requires two evaluations of the
function and one evaluation of its first derivative per iteration
and it has the efficiency index equal to $3^{\frac{1}{3}}\approx
1.44225$. We describe the analysis of the proposed methods along
with numerical experiments including comparison with existing
methods. Moreover, the attraction basins are shown and compared
with other existing methods.

\textbf{Keywords}: Multi-point iterative methods, Potra-Pt{\'a}k
method,  Multiple roots, Basin of attraction.
\end{abstract}

\section{Introduction }
\label{intro} Solving nonlinear equations based on iterative
methods is a basic and extremely valuable tool in all fields of
science as well as economics and engineering. The important
aspects related to these methods are order of convergence and
number of function evaluations. Therefore, it is favorable to
attain the highest possible convergence order with fixed number of
function evaluations for each iteration. The aim of the paper is
to modify third-order the Potra-Pt{\'a}k method to solve nonlinear
equations for multiple zeros with same order of convergence and
efficiency index. The efficiency index of an iterative method of
order $p$ requiring $k$ function evaluations per iteration is
defined by $E(k,p)=\sqrt[k]{p}$, see~\cite{Ostrowski}.

Let $\alpha$ be multi roots of $f(x)=0$ with multiplicity $m$
i.e., $f^{(i)}(\alpha)=0, \quad i=0,1,\ldots,m-1$ and
$f^{(m)}(\alpha) \neq 0$. If functions $f^{(m-1)}$ and $f^{1/m}$
have only a simple zero at $\alpha$, any of the iterative methods
for a simple zero may be used \cite{chun,Traub}. The modified
Newton method defined in \cite{scheroder} is one of the most well
known iterative methods for multiple roots
\begin{equation*}
x_{n+1}=x_n-m\frac{f(x_n)}{f'(x_n)},
\end{equation*}
which converges quadratically.

In the recent years, a wide collection of iterative methods for
finding simple roots and multiple roots of nonlinear equations
have been presented in some journals, see
\cite{Lotfi1,Petkovic,Sharifi,Sharifi2}. In order to improve the
convergence of iterative methods for multiple roots, many
researchers such as Chun et al. \cite{chun}, Dong \cite{Dong},
Hansen and Patrick \cite{Halley}, Heydari et al. \cite{Heydari},
Osada \cite{Osada}, Victory and Neta \cite{victor} proposed
various iterative methods when the multiplicity $m$ is known.

This paper is organized as follow: Section \ref{sec:2} is devoted
to the construction and convergence analysis of a new method which
theoretical proof has been given to reveal the third-order of
convergence. In Section \ref{sec:3}, different numerical tests
confirm the theoretical results and allow us to compare this
method with other known methods. Comparisons of attraction basins
with other methods are illustrated in this section as well.
Finally, a conclusion is provided in Section \ref{sec:conclusion}.


\section{The method and analysis of convergence} \label{sec:2}

\subsection{Extension for multiple roots}
In this section, we proceed to develop a new iterative method to
find an approximate root for a nonlinear
equation $f(x)=0$ where $f:\mathbb{C}\rightarrow\mathbb{C}$.\\
Our main aim is to extend Potra-Pt{\'a}k iterative method for
multiple roots and build an iterative method using some parameters
 without any additional evaluations
of the function or its derivatives.\\
The Potra-Pt{\'a}k iterative method \cite{Potra} is given by
\begin{equation}
  \begin{array}{lrl}\label{c1}
 y_{n}=x_{n}-\frac{f(x_{n})}{f'(x_{n})},&\\
  [1ex]
x_{n+1}=x_n-\frac{f(x_n)+f(y_n)}{f'(x_n)},\quad(n=0,1,\ldots),
      \end{array}
\end{equation}
with the initial approximation of $x_0$ sufficiently close to
$x^*$. Convergence order of Potra-Pt{\'a}k method for
approximating simple zero of nonlinear equations is three, whereas
for finding multiple zeros is linear. In order to do so, we add
two parameters $\alpha$ and $\beta$ on the second term of
(\ref{c1}). So
\begin{equation}
  \begin{array}{lrl}\label{c2}
 y_{n}=x_{n}-\frac{f(x_{n})}{f'(x_{n})},&\\
  [1ex]
x_{n+1}=x_n-\frac{\alpha f(x_n)+f(y_n)}{\beta f'(x_n)},
      \end{array}
\end{equation}
where $\alpha$ and $\beta$ are parameters to be chosen for maximal order of convergence.
\subsection{Proof of convergence}
In what follow, we describe the convergence analysis on iterative
method (\ref{c2}).\\
\begin{thm}\label{111}
Let $x^* \in D$  be a multiple zero of a sufficiently
differentiable function $f:D \subset \mathbb{C} \rightarrow
\mathbb{C}$ for an open interval $D$ with integer multiplicity
$m\geq1$, which includes $x_0$ as an initial approximation of
$x^*$. Then, method (\ref{c2}) has order of three,
$\alpha=\frac{(\mu-1)\mu^{m-1}}{m^m}$ and
$\beta=\frac{\mu^{m-1}}{m^{m+1}}$ where $\mu=m-1$.\\
\end{thm}
\begin{proof}
Let $e_{n}:=x_{n}-\alpha$, $e_{n,y}:=y_{n}-\alpha$,
$c_{i}:=\frac{m!}{(m+i)!}\frac{f^{(m+i)}(\alpha)}{f^{(m)}(\alpha)}$.
Using the fact that $f(\alpha)=0$, Taylor expansion of $f$ at
$\alpha$ yields
\begin{equation}\label{a1}
f(x_{n}) =e_n^m\left(c_{0} +c_1e_n+
c_{2}e_{n}^{2}+c_{3}e_{n}^{3}\right)+O(e_n^{4}),
\end{equation}
and
\begin{equation}\label{a2}
f^{'}(x_{n})
=e_n^{m-1}(m+(m+1)c_1e_n+(m+2)c_2e_n^2+(m+3)c_3e_n^3+O(e_n^4)).
\end{equation}
Hence
\begin{equation}\label{a4}
e_{n,y}=y_{n}-\alpha=
\frac{-1+m}{m}e_n-\frac{c_1}{m^2c_0}e_n^2+\frac{-(1+m)c_1^2+2mc_0c_2}{m^3c_0^2}e_n^3+O(e_n^{4}).
\end{equation}
For $f(y_n)$ we have
\begin{equation}\label{a5}
f(y_{n}) =e_{n,y}^m\left(c_{0} +c_1e_{n,y}+
c_{2}e_{n,y}^{2}+c_{3}e_{n,y}^{3}\right)+O(e_{n,y}^{4}).
\end{equation}
Substituting (\ref{a1})-(\ref{a5}) in (\ref{c2}), we obtain
\begin{equation*}
e_{n+1}= D_1e_n+D_2e_n^2+D_3e_n^3+O(e_n^4),
\end{equation*}
where
\begin{equation*}\label{a6}
D_1=1-\frac{\alpha+(\mu/m)^3}{\beta m},\quad D_2=\frac{(\mu
\alpha+(\mu-1)(\mu/m)^m)c_1}{\beta \mu m^2 c_0},
\end{equation*}
and
\begin{equation*}\label{a8}
D_3=\frac{-2c_1^2 \alpha (-1+m^2)+4\alpha \mu c_0 c_2+c_1^2(\mu/m)^m ((6+(1-2\mu)m)+2c_0c_2(1+m(2\mu-3)))}{2\beta \mu m^3 c_0^2}.
\end{equation*}
Therefor, to provide the order of convergence three, it is
necessary to choose $D_i=0 \quad (i=1,2)$, which gives
\begin{equation*}\label{a8}
\alpha=(\mu-1)\frac{\mu^{\mu}}{m^m},\quad  \text{and} \quad \beta=\frac{\mu^{\mu}}{m^{m+1}},
\end{equation*}
where $\mu=m-1$ and the error equation becomes
\begin{equation*}
e_{n+1}=\frac{(2+m)c_1^2-2\mu c_0c_2}{2m^2c_0^2}e_n^3+O(e_n^4).
\end{equation*}
\end{proof}
Therefore method (\ref{c2}) has convergence order three, which
proved here.

\section{Results and discussions}\label{sec:3}
 In this section we check the effectiveness of the modified
 Potra-Pt{\'a}k method and compare with other existing methods
 which have the same order of convergence.\\
\subsection{Numerical results}
Substituting $\alpha$ and $\beta$ in the Theorem \ref{111} into
(2.2), the modified Potra-Pt{\'a}k is
\begin{equation}
  \begin{array}{lrl}\label{c4}
 y_{n}=x_{n}-\frac{f(x_{n})}{f'(x_{n})},&\\
  [1ex]
x_{n+1}=x_n+m\frac{(\mu-1)\mu^{\mu}f(x_n)-m^mf(y_n)}{\mu^{\mu}f'(x_n)},
      \end{array}
\end{equation}
where $\mu=m-1$.\\
The Osada's method \cite{Osada}, is given by
\begin{equation}\label{m1}
x_{n+1}=x_n-\frac{1}{2}m(m+1)\frac{f(x_n)}{f'(x_n)}+\frac{1}{2}(m-1)^2\frac{f'(x_n)}{f''(x_n)}.
\end{equation}
The Dong's method \cite{Dong}, is given by
\begin{equation}
  \begin{array}{lrl}\label{m2}
 y_n=x_n+\sqrt{m}\frac{f(x_n)}{f'(x_n)},&\\
  [1ex]
x_{n+1}=y_n-m(1-\frac{1}{\sqrt{m}})^{1-m}\frac{f(y_n)}{f'(x_n)}.
  \end{array}
\end{equation}
The Chun's method \cite{chun}, is given by
\begin{equation}
  \begin{array}{lrl}\label{m3}
 x_{n+1}=x_n-\frac{m((2\gamma-1)m+3-2\gamma)}{2}\frac{f(x_n)}{f'(x_n)}&\\
  [1ex]
~~~~~~+\frac{\gamma(m-1)^2}{2}\frac{f'(x_n)}{f''(x_n)}-\frac{(1-\gamma)m^2}{2}\frac{f(x_n)^2f''(x_n)}{f'(x_n)^3}.
  \end{array}
\end{equation}
In the numerical experiments of this paper we use
$\gamma=-1$.\\
We test method (\ref{c2}) on a number of nonlinear equations.
Numerical computations have been carried out using variable
precision arithmetic with 100 significant decimal digits in the
programming package of Mathematica 8 \cite{Hazrat}.
\begin{table}[!ht]
\begin{center}
\begin{tabular}{c l l l l}
\hline
\\
 \multicolumn{5}{l}{$ f_1(x)=\big(\ln(1+x^2)+\mathrm{e}^{x^2-3x}\sin x\big)^6,\quad m=6\quad x^*=0,\quad x_0=0.3$} \\
 \hline
&Method (\ref{c4})&Method (\ref{m1})&Method (\ref{m2})&Method (\ref{m3})\\
 \hline
$\vert x_{1}-x^{*} \vert$ & $0.656\mathrm{e}{-}1$ & $0.678\mathrm{e}{-}1$ & $0.552\mathrm{e}{-}1$ & $0.736\mathrm{e}{-}1$  \\
$\vert x_{2}-x^{*} \vert$ & $0.128\mathrm{e}{-}2$& $0.178\mathrm{e}{-}2$ & $0.511\mathrm{e}{-}3$ & $0.118\mathrm{e}{-}2$  \\
$\vert x_{3}-x^{*} \vert$ & $0.696\mathrm{e}{-}8$ & $0.247\mathrm{e}{-}7$ & $0.273\mathrm{e}{-}9$ & $0.186\mathrm{e}{-}8$ \\
COC & $3.0000$ & $3.0000$ & $3.0002$ & $3.0020$\\
ACOC & $3.1010$ & $3.0997$ & $3.0914$ & $3.2527$\\
\hline
\\
 \multicolumn{5}{l}{$ f_2(x)=\big(x^3+\ln(1+x)\big)^7,\quad m=7\quad x^*=0,\quad x_0=0.2$} \\
 \hline
 &Method (\ref{c4})&Method (\ref{m1})&Method (\ref{m2})&Method (\ref{m3})\\
 \hline
$\vert x_{1}-x^{*} \vert$ & $0.880\mathrm{e}{-}2$ & $0.792\mathrm{e}{-}4$ & $0.122\mathrm{e}{-}4$ & $0.987\mathrm{e}{-}4$  \\
$\vert x_{2}-x^{*} \vert$ & $0.417\mathrm{e}{-}6$& $0.139\mathrm{e}{-}15$ & $0.742\mathrm{e}{-}19$ & $0.340\mathrm{e}{-}15$  \\
$\vert x_{3}-x^{*} \vert$ & $0.481\mathrm{e}{-}19$ & $0.918\mathrm{e}{-}39$ & $0.918\mathrm{e}{-}39$ & $0.918\mathrm{e}{-}39$ \\
COC & $3.0000$ & $3.0000$ & $3.0000$ & $3.0000$\\
ACOC & $2.9921$ & $3.0000$ & $3.0000$ & $3.0000$\\
\hline
\\
 \multicolumn{5}{l}{$ f_3(x)=(x^6-8)^2\ln(x^6-7),\quad m=3\quad x^*=\sqrt{2},\quad x_0=1.5$} \\
 \hline
 &Method (\ref{c4})&Method (\ref{m1})&Method (\ref{m2})&Method (\ref{m3})\\
 \hline
$\vert x_{1}-x^{*} \vert$ & $0.328\mathrm{e}{-}2$ & $0.414\mathrm{e}{-}2$ & $0.236\mathrm{e}{-}2$ & $0.422\mathrm{e}{-}2$  \\
$\vert x_{2}-x^{*} \vert$ & $0.914\mathrm{e}{-}6$& $0.247\mathrm{e}{-}5$ & $0.200\mathrm{e}{-}6$ & $0.591\mathrm{e}{-}5$  \\
$\vert x_{3}-x^{*} \vert$ & $0.184\mathrm{e}{-}17$ & $0.491\mathrm{e}{-}15$ & $0.117\mathrm{e}{-}16$ & $0.129\mathrm{e}{-}16$ \\
COC & $3.0000$ & $3.0000$ & $3.0000$ & $3.0000$\\
ACOC & $3.0085$ & $3.0067$ & $3.0051$ & $3.0339$\\
\hline
\\
 \multicolumn{5}{l}{$ f_4(x)=\big(\ln(x^2-x+1)+4\sin (x-1)\big)^{10},\quad m=10\quad x^*=1,\quad x_0=1.2$} \\
 \hline
 &Method (\ref{c4})&Method (\ref{m1})&Method (\ref{m2})&Method (\ref{m3})\\
 \hline
$\vert x_{1}-x^{*} \vert$ & $0.102\mathrm{e}{-}2$ & $0.182\mathrm{e}{-}2$ & $0.126\mathrm{e}{-}2$ & $0.102\mathrm{e}{-}2$  \\
$\vert x_{2}-x^{*} \vert$ & $0.583\mathrm{e}{-}10$& $0.175\mathrm{e}{-}8$ & $0.393\mathrm{e}{-}9$ & $0.583\mathrm{e}{-}10$  \\
$\vert x_{3}-x^{*} \vert$ & $0.107\mathrm{e}{-}31$ & $0.155\mathrm{e}{-}26$ & $0.119\mathrm{e}{-}29$ & $0.107\mathrm{e}{-}31$ \\
COC & $3.0000$ & $3.0000$ & $3.0000$ & $3.0000$\\
ACOC & $3.0003$ & $2.9998$ & $2.9999$ & $3.0003$\\
\hline
\\
 \multicolumn{5}{l}{$ f_5(x)=\ln^2(x-2)(\mathrm{e}^{x-3}-1)\sin\big(\pi x/3\big),\quad m=4\quad x^*=3,\quad x_0=3.1$} \\
 \hline
 &Method (\ref{c4})&Method (\ref{m1})&Method (\ref{m2})&Method (\ref{m3})\\
 \hline
$\vert x_{1}-x^{*} \vert$ & $0.126\mathrm{e}{-}3$ & $0.144\mathrm{e}{-}3$ & $0.890\mathrm{e}{-}4$ & $0.253\mathrm{e}{-}3$  \\
$\vert x_{2}-x^{*} \vert$ & $0.572\mathrm{e}{-}13$& $0.106\mathrm{e}{-}12$ & $0.132\mathrm{e}{-}13$ & $0.905\mathrm{e}{-}12$  \\
$\vert x_{3}-x^{*} \vert$ & $0.531\mathrm{e}{-}41$ & $0.424\mathrm{e}{-}40$ & $0.493\mathrm{e}{-}43$ & $0.414\mathrm{e}{-}37$ \\
COC & $3.0000$ & $3.0000$ & $3.0000$ & $3.0000$\\
ACOC & $3.0003$ & $3.0000$ & $3.0000$ & $3.0000$\\
\hline
\end{tabular}
\end{center}
\vspace*{-3ex} \caption{Errors, COC and ACOC for methods
(\ref{c4}), (\ref{m1}), (\ref{m2}) and (\ref{m3}). \label{table2}}
\end{table}
\newpage
In Table \ref{table2}, new modified Potra-Pt{\'a}k method
(\ref{c4}) is compared with the methods (\ref{m1}),(\ref{m2}) and
(\ref{m3}) on five nonlinear equations which are illustrated from
left to right respectively. In addition the computational order of
convergence is defined by \cite{coc}
\begin{equation*}
\textup{COC}\approx\frac{\ln|(x_{n+1}-\alpha)/(x_{n}-\alpha)|}{\ln|(x_{n}-\alpha)/(x_{n-1}-\alpha)|},
\end{equation*}
 and the approximated computational order of convergence
 is defined by \cite{acoc}
\begin{equation*}
\textup{ACOC}\approx\frac{\ln|(x_{n+1}-x_{n})/(x_{n}-x_{n-1})|}{\ln|(x_{n}-x_{n-1})/(x_{n-1}-x_{n-2})|}.
\end{equation*}
Table \ref{table2} show that new modified Potra-Pt{\'a}k method
\ref{c4} support the established theorem given in the previous
section.

\subsection{Comparison of attraction basins}
In this section, we check the comparison of iterative methods in
the complex plane by using basins of attraction. The basin of
attraction is a method to visually comprehend how an algorithm
behaves as a function of the various starting points. With
attractions basins, the study of dynamic behavior of the rational
functions associated to an iterative method gives important
information about the convergence and stability of the scheme
\cite{Lotfi1,Sharifi}. Let $G:\mathbb{C} \to \mathbb{C} $ be a
rational mapping on the complex plane. For $z\in \mathbb{C} $, we
define its orbit as the set $orb(z)=\{z,\,G(z),\,G^2(z),\dots\}$.
A point $z_0 \in \mathbb{C} $ is called periodic point with
minimal period $m$ if $G^m(z_0)=z_0$, where $m$ is the smallest
integer with this property. A periodic point with minimal period
$1$ is called fixed point. Moreover, a point $z_0$ is called
attracting if $|G'(z_0)|<1$, repelling if $|G'(z_0)|>1$, and
neutral otherwise. The Julia set of a nonlinear map $G(z)$,
denoted by $J(G)$, is the closure of the set of its repelling
periodic points. The complement of $J(G)$ is the Fatou set $F(G)$,
where the basin of attraction of the different roots lie
\cite{Babajee}.

We use the basin of attraction for comparing the iteration
algorithms. For the dynamical point of view, we take a $256 \times
256$ grid of the square $[-3,3]\times[-3,3]\in \mathbb{C}$ and
assign a color to each point $z_0\in D$ according to the simple
root to which the corresponding orbit of the iterative method
starting from $z_0$ converges, and we mark the point as black if
the orbit does not converge to a root, in the sense that after at
most 100 iterations it has a distance to any of the roots, which
is larger than $10^{-3}$. In this way, we distinguish the
attraction basins by their color for different methods.

In the following figures, the roots of each function are assigned
to different colors. In the basins of attraction, the number of
iterations needed to approach the solution is indicated in darker
or brighter colors. Black color denotes lack of convergence to any
of the roots or convergence to infinity.

\begin{figure}[ht!]
\begin{minipage}[b]{0.28\linewidth}
\centering
\includegraphics[width=\textwidth]{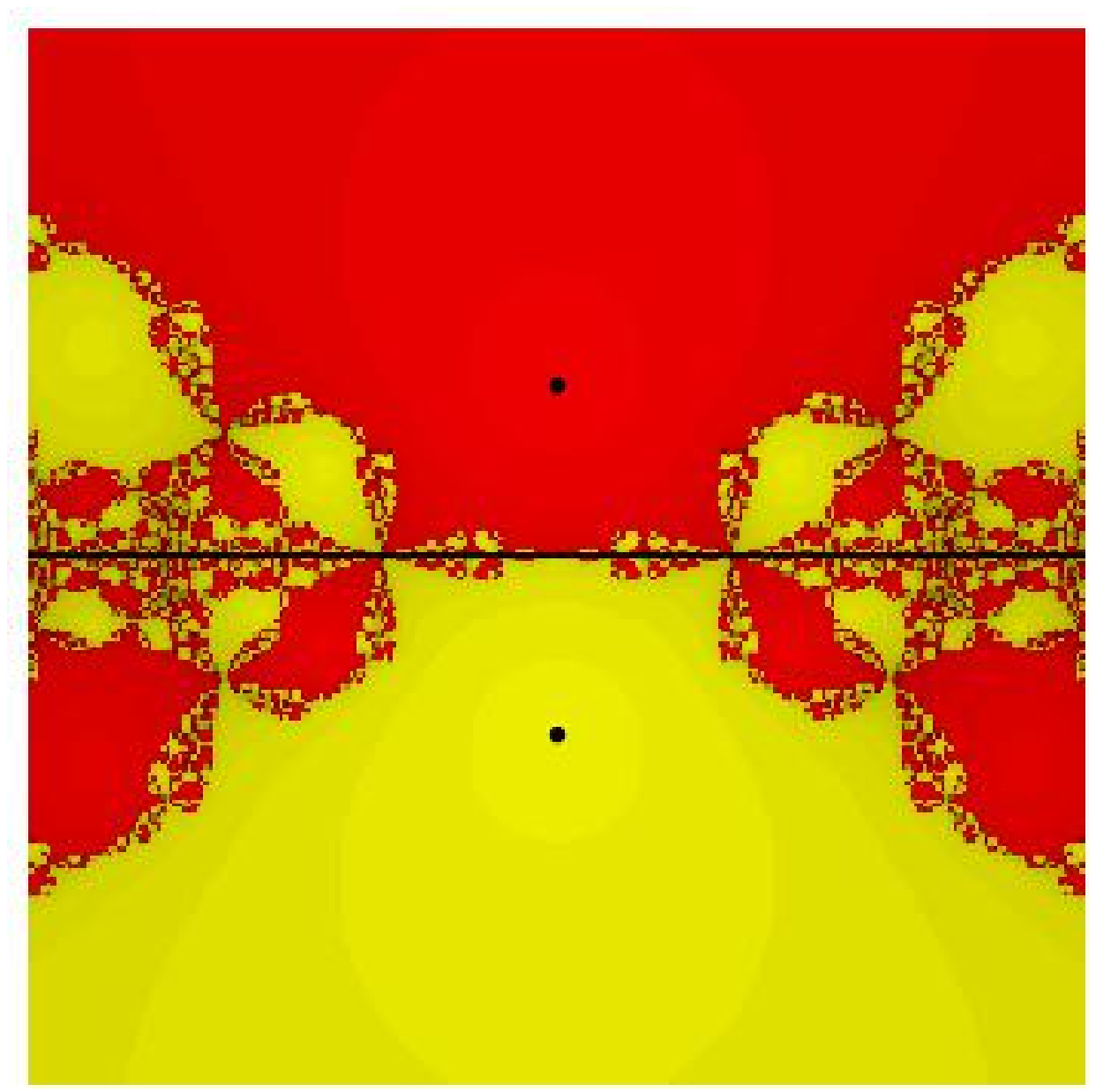}
\end{minipage}
\hspace{0.5cm}
\begin{minipage}[b]{0.28\linewidth}
\centering
\includegraphics[width=\textwidth]{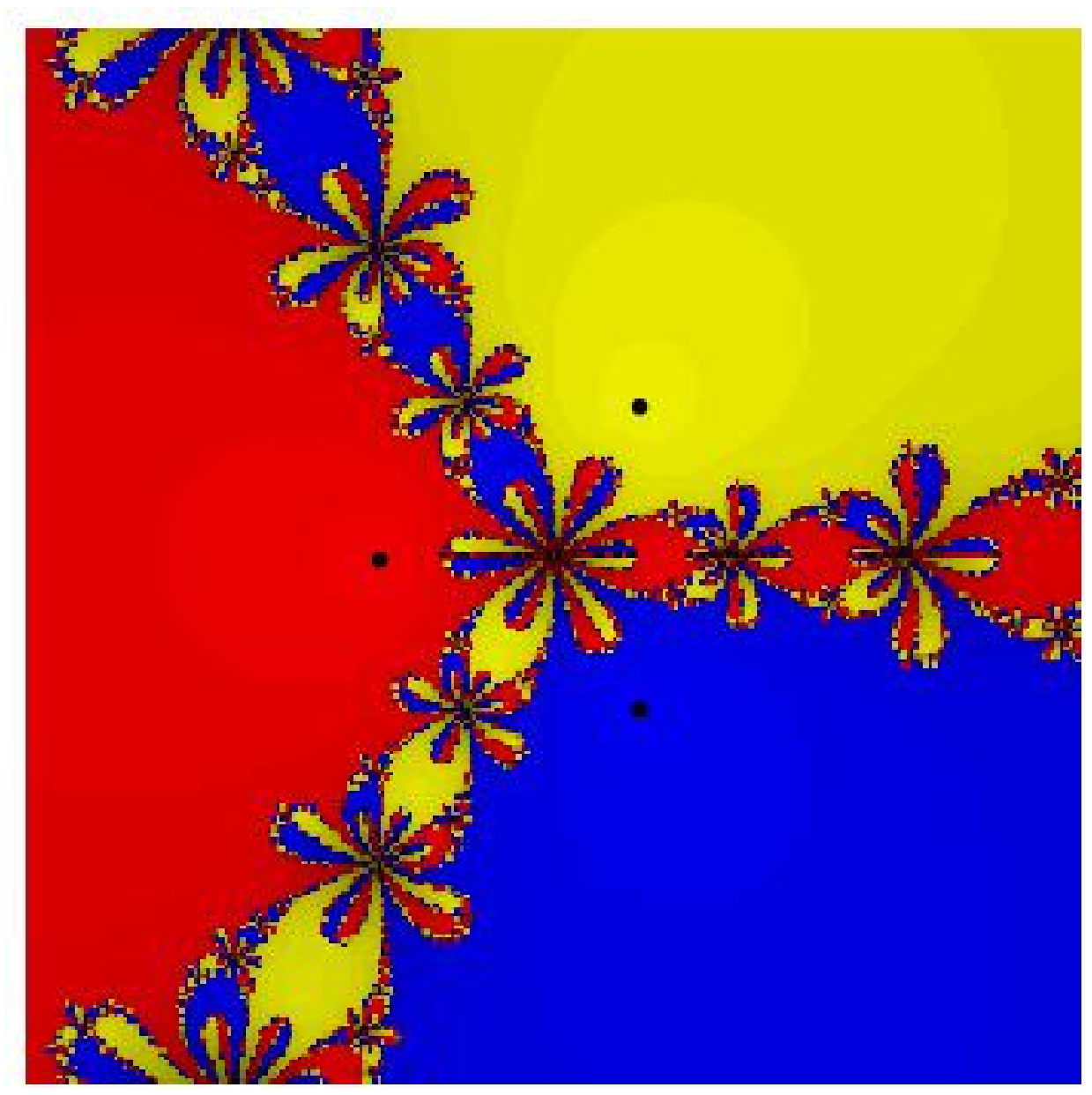}
\end{minipage}
\hspace{0.5cm}
\begin{minipage}[b]{0.28\linewidth}
\centering
\includegraphics[width=\textwidth]{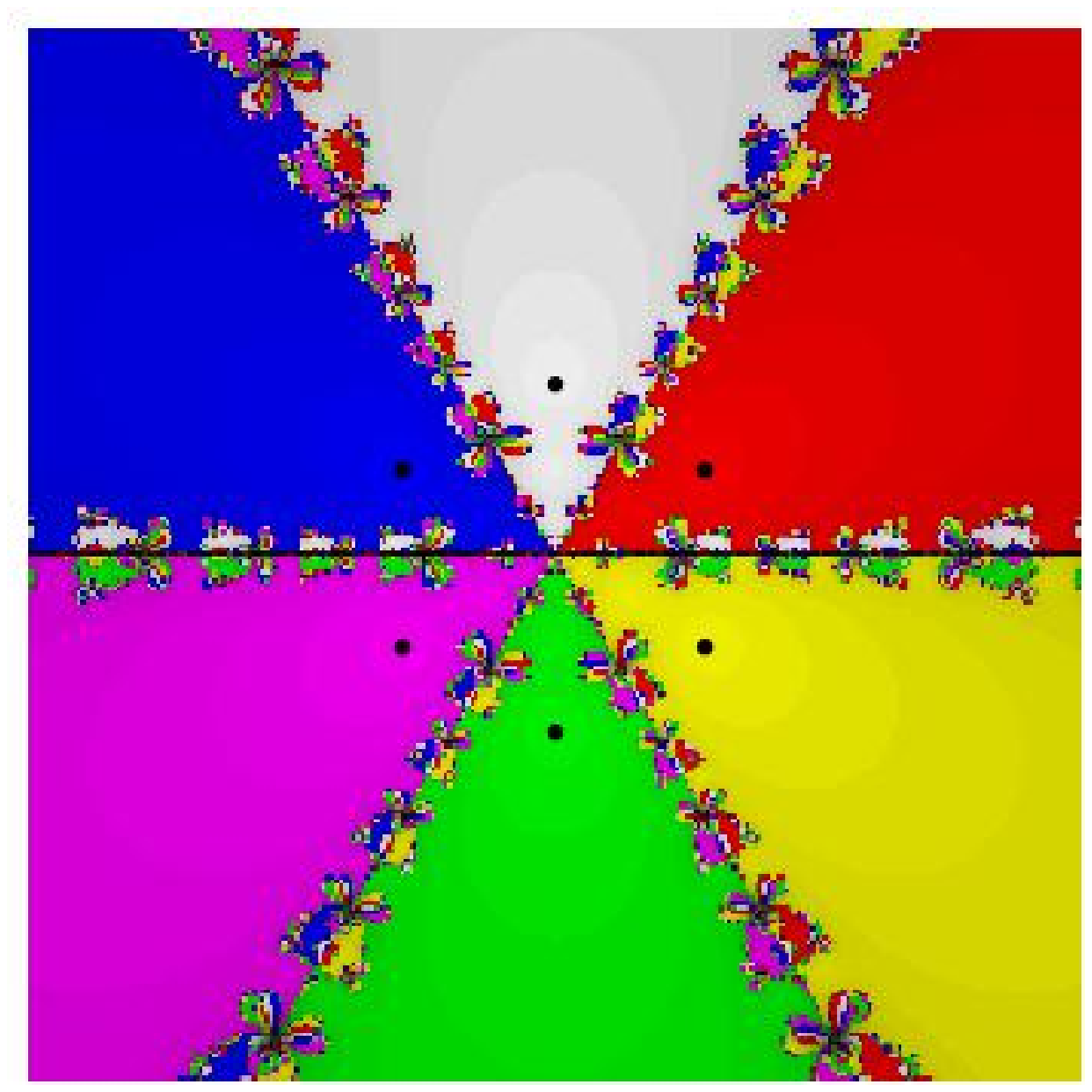}
\end{minipage}
\caption{Method (\ref{c1}) for test problems
$p_1(z)=z+\frac{1}{z}$, $p_2(z)=z^3+1$ and
$p_3(z)=z^5+\frac{1}{z}$ respectively, without any multiplicity}
\label{fig:figure1}
\end{figure}

\begin{figure}[ht!]
\begin{minipage}[b]{0.28\linewidth}
\centering
\includegraphics[width=\textwidth]{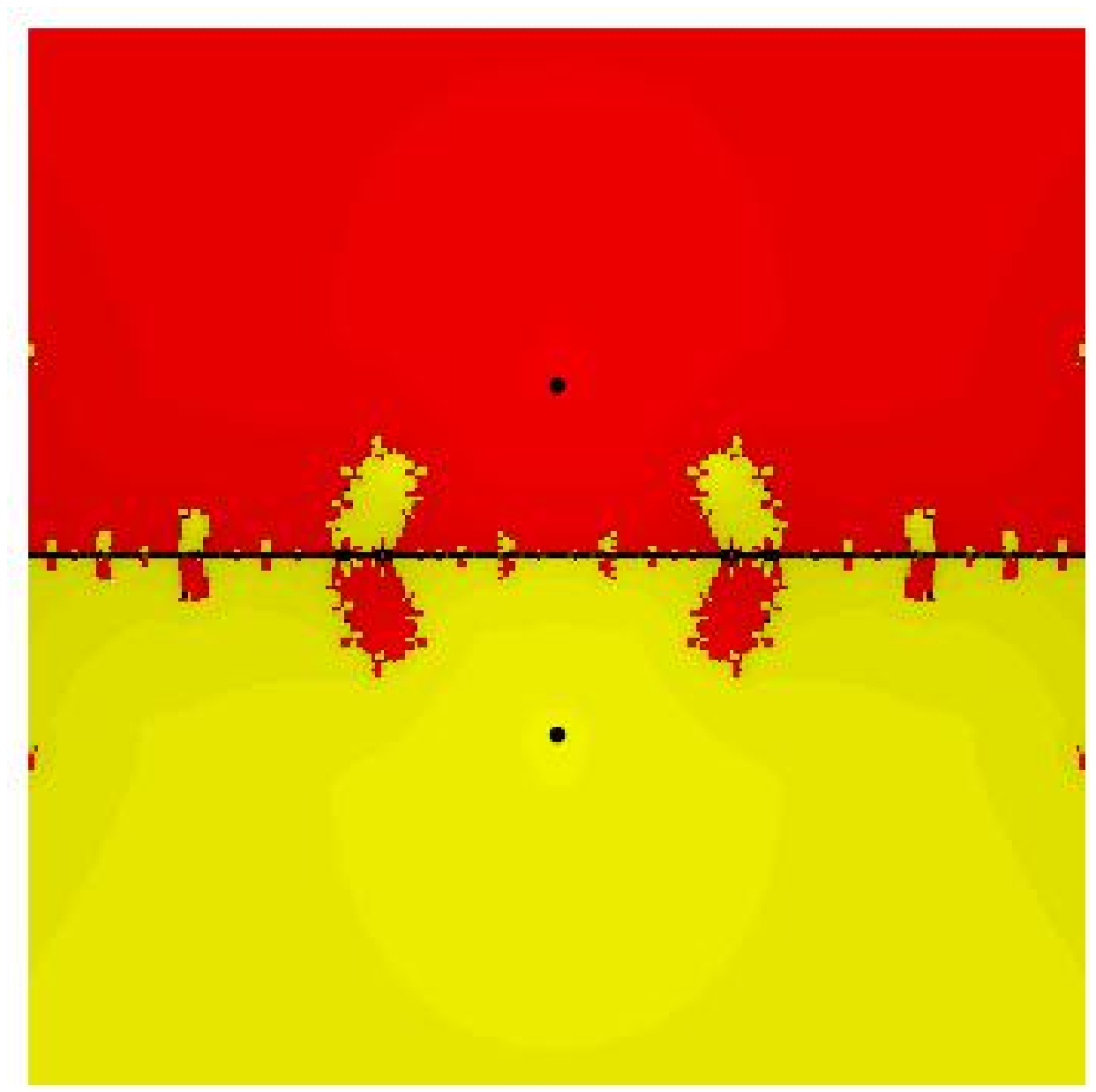}
\end{minipage}
\hspace{0.5cm}
\begin{minipage}[b]{0.28\linewidth}
\centering
\includegraphics[width=\textwidth]{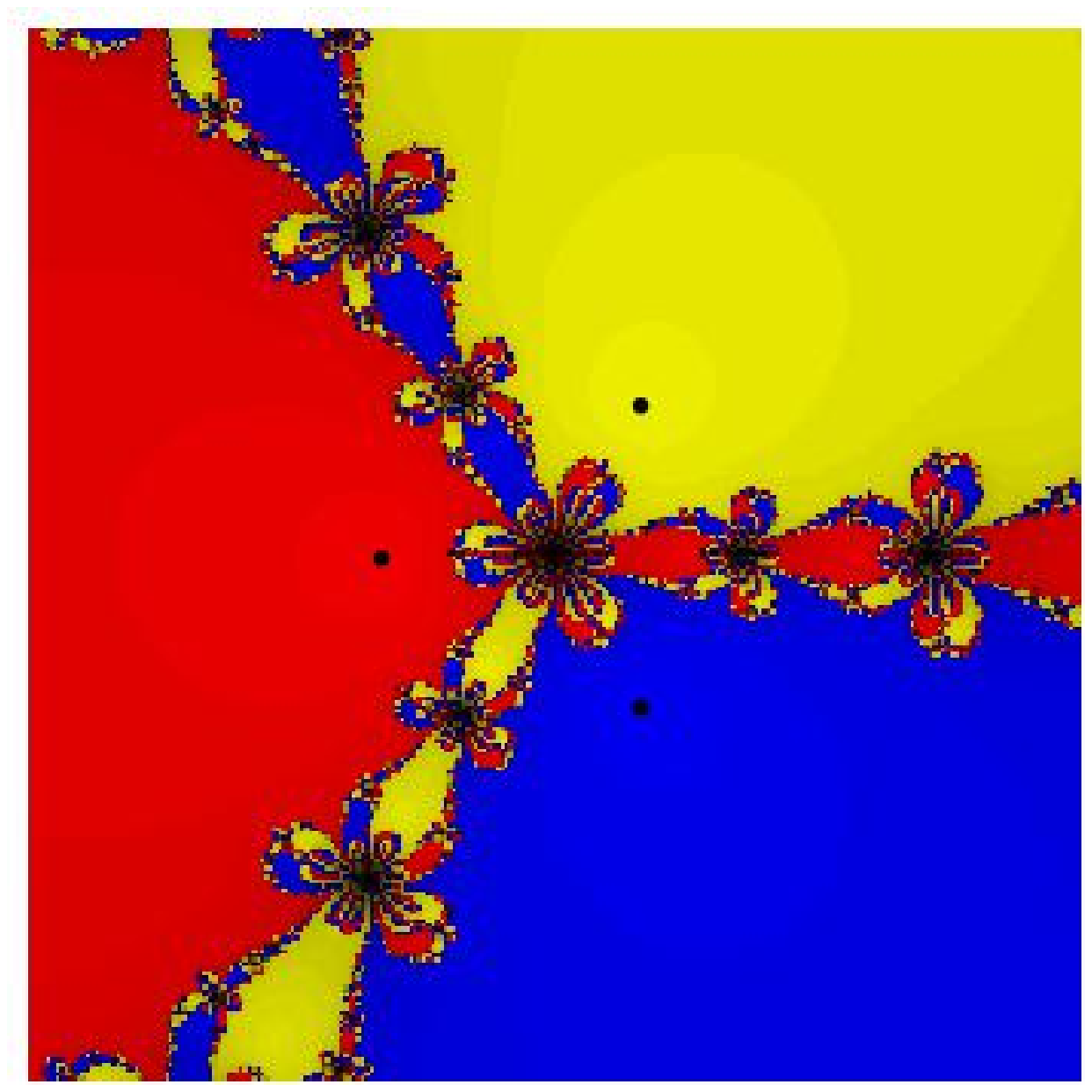}
\end{minipage}
\hspace{0.5cm}
\begin{minipage}[b]{0.28\linewidth}
\centering
\includegraphics[width=\textwidth]{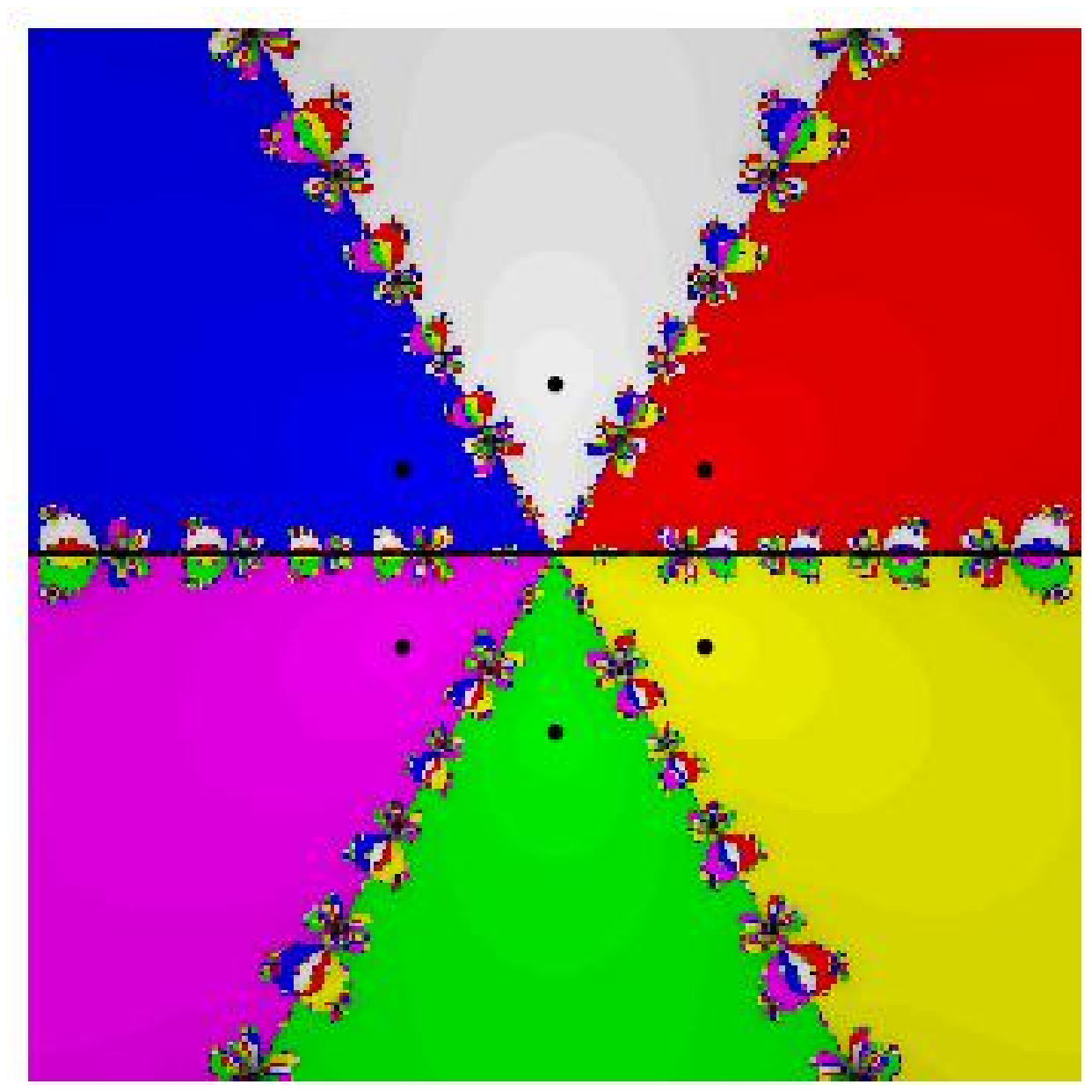}
\end{minipage}
\caption{Method (\ref{c4}) for test problems
$p_1(z)=\left(z+\frac{1}{z}\right)^5$, $p_2(z)=(z^3+1)^3$ and
$p_3(z)=\left(z+\frac{1}{z}\right)^2$ respectively}
\label{fig:figure2}
\end{figure}

\begin{figure}[ht!]
\begin{minipage}[b]{0.28\linewidth}
\centering
\includegraphics[width=\textwidth]{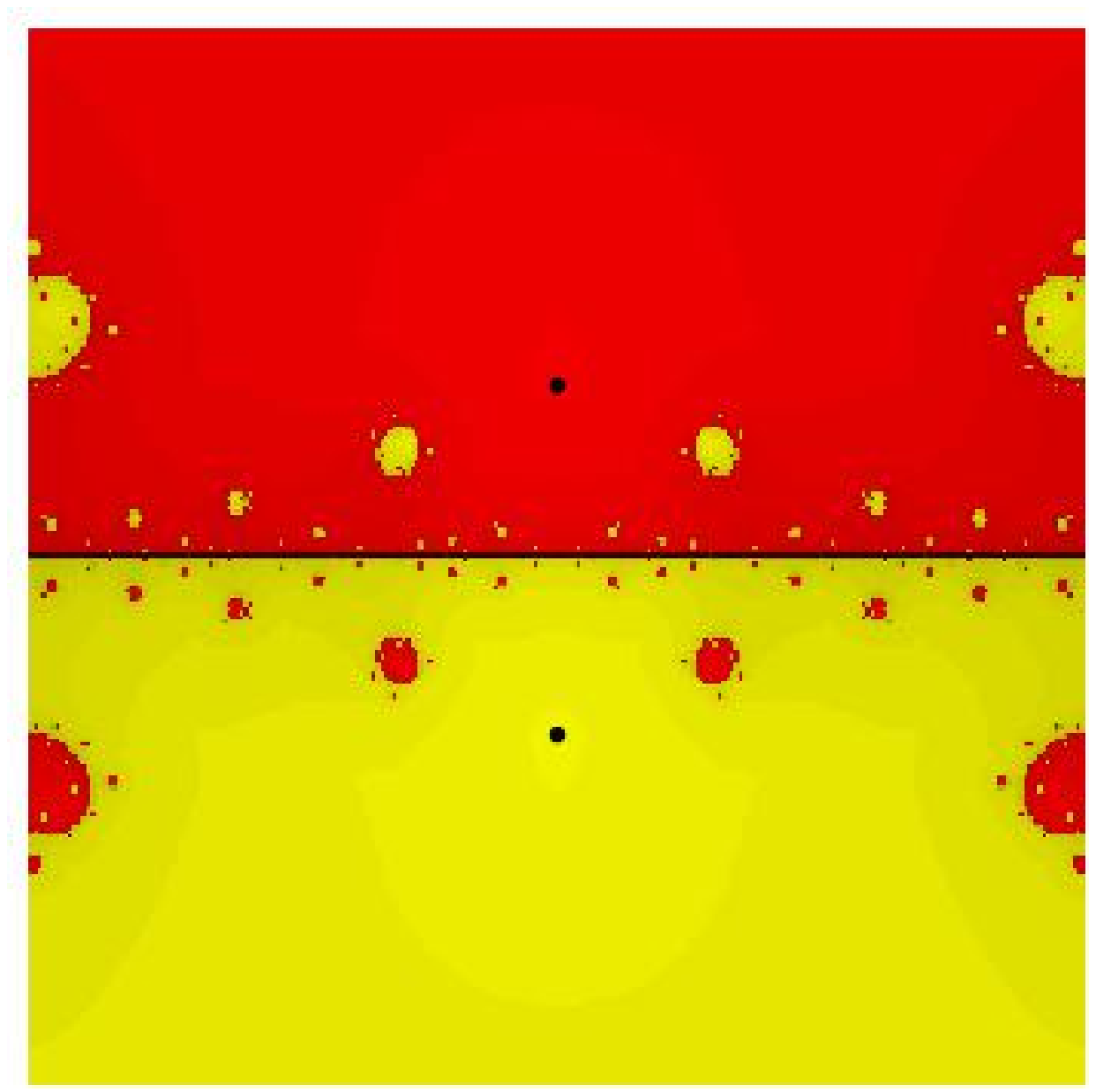}
\end{minipage}
\hspace{0.5cm}
\begin{minipage}[b]{0.28\linewidth}
\centering
\includegraphics[width=\textwidth]{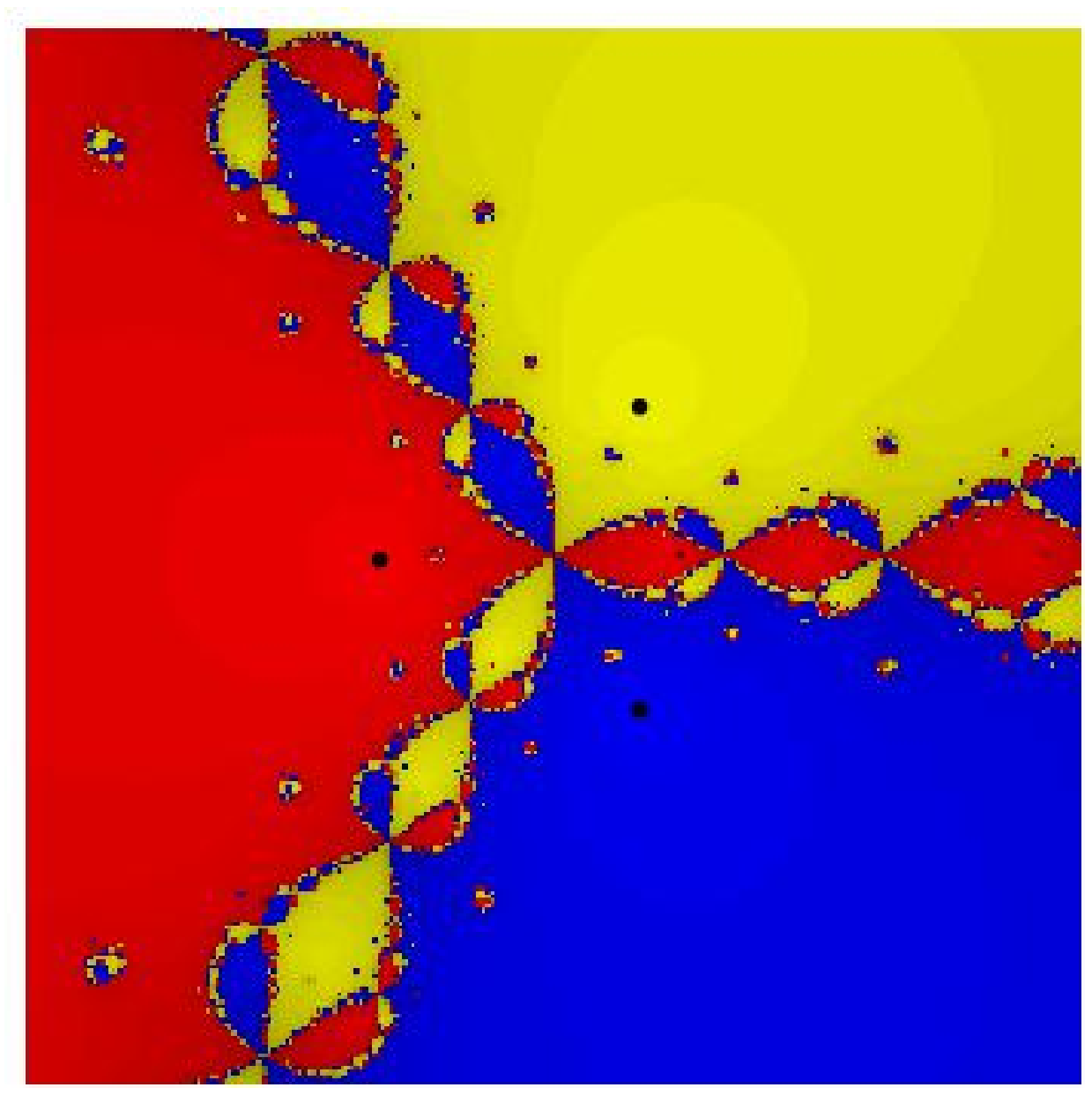}
\end{minipage}
\hspace{0.5cm}
\begin{minipage}[b]{0.28\linewidth}
\centering
\includegraphics[width=\textwidth]{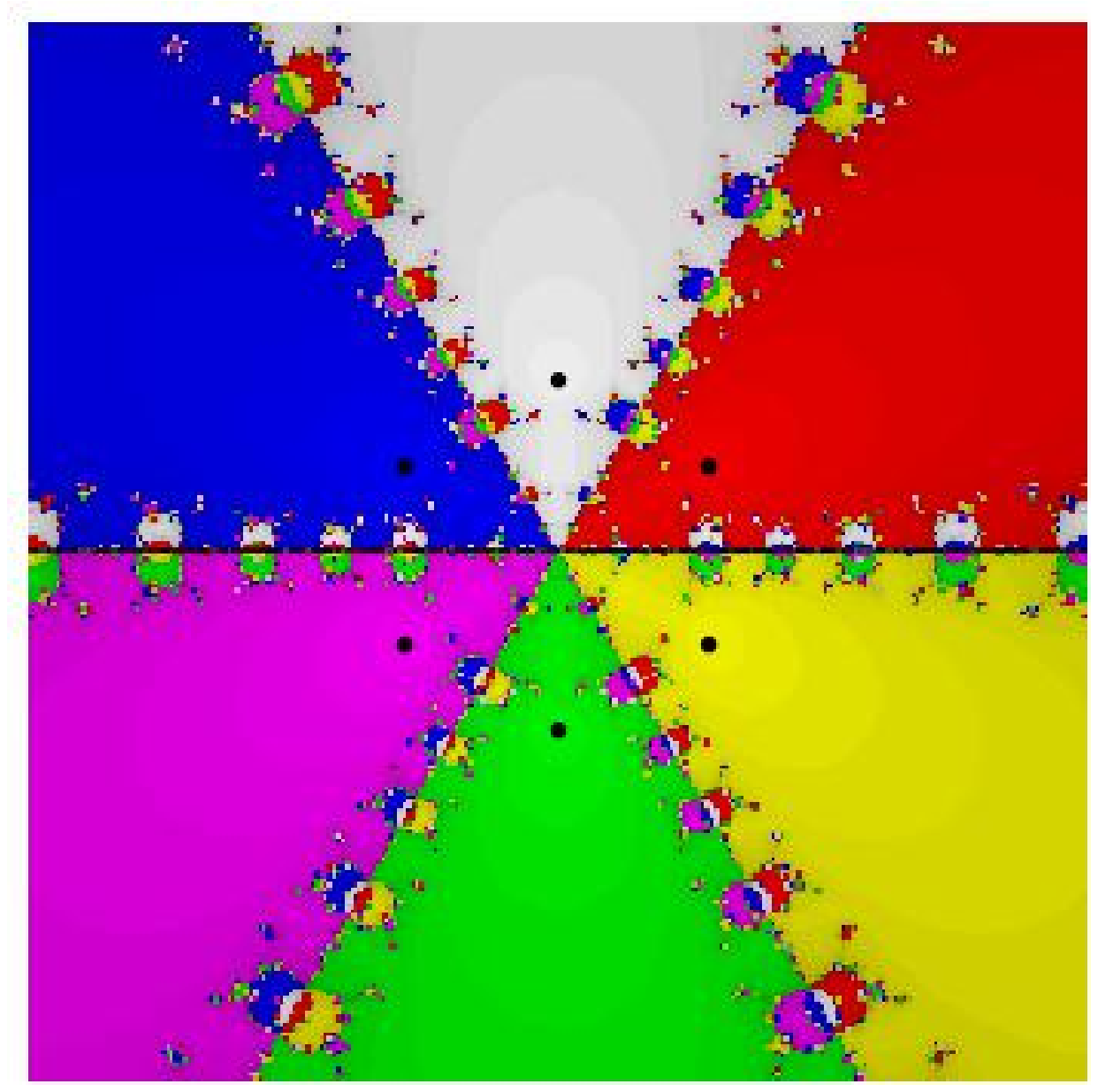}
\end{minipage}
\caption{Method (\ref{m1}) for test problems
$p_1(z)=\left(z+\frac{1}{z}\right)^5$, $p_2(z)=(z^3+1)^3$ and
$p_3(z)=\left(z+\frac{1}{z}\right)^2$ respectively}
\label{fig:figure3}
\end{figure}

\begin{figure}[ht!]
\begin{minipage}[b]{0.28\linewidth}
\centering
\includegraphics[width=\textwidth]{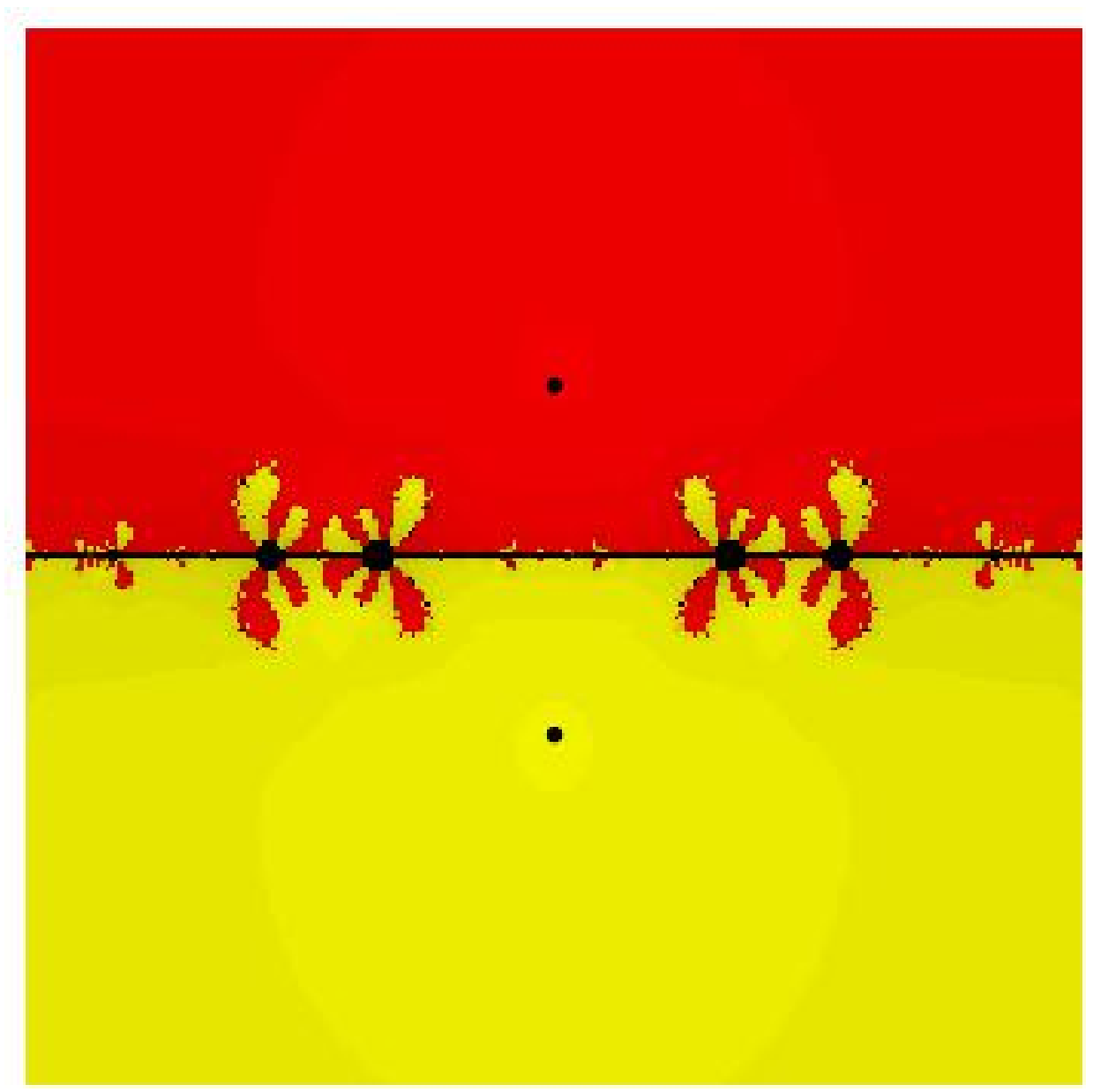}
\end{minipage}
\hspace{0.5cm}
\begin{minipage}[b]{0.28\linewidth}
\centering
\includegraphics[width=\textwidth]{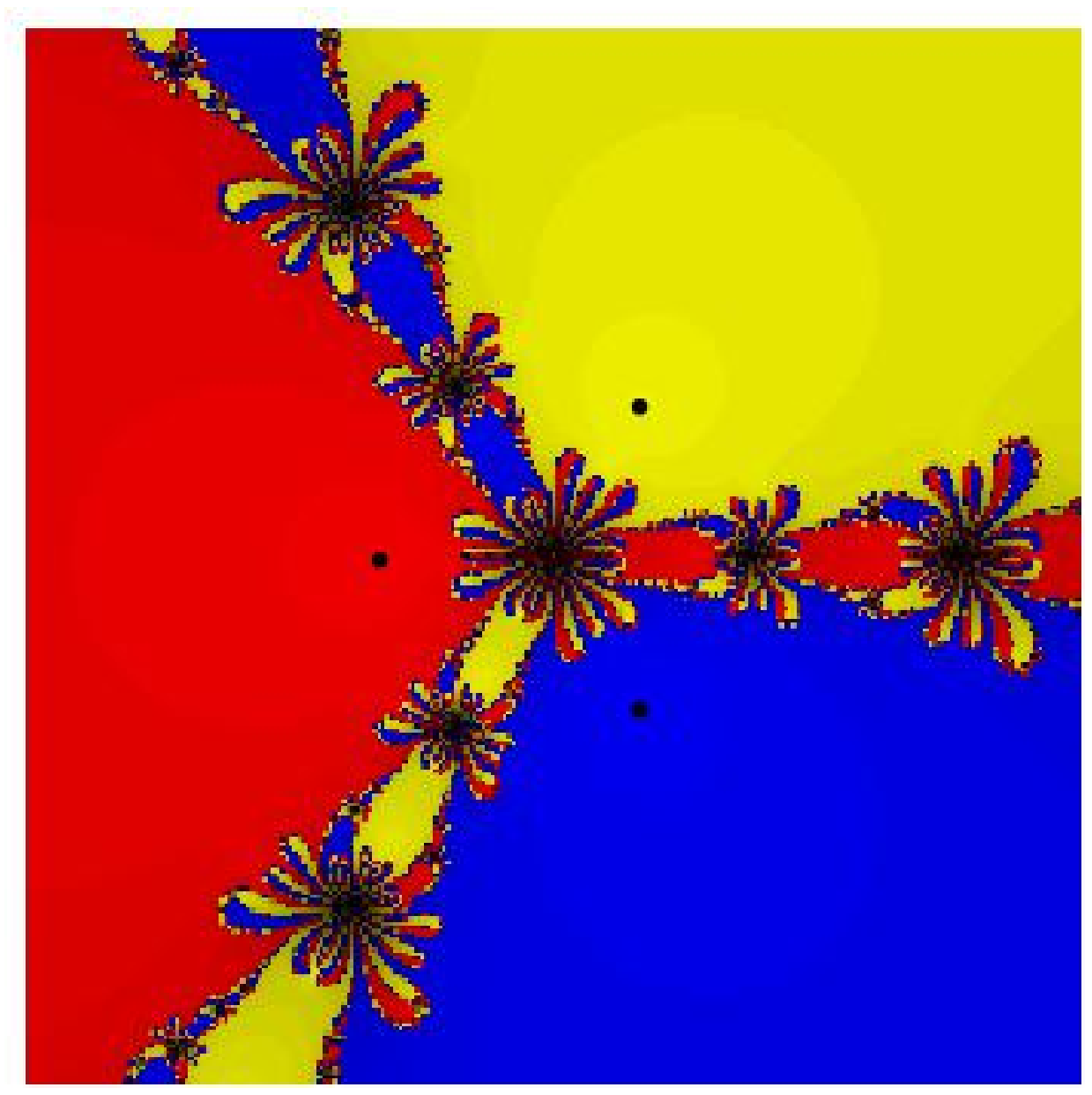}
\end{minipage}
\hspace{0.5cm}
\begin{minipage}[b]{0.28\linewidth}
\centering
\includegraphics[width=\textwidth]{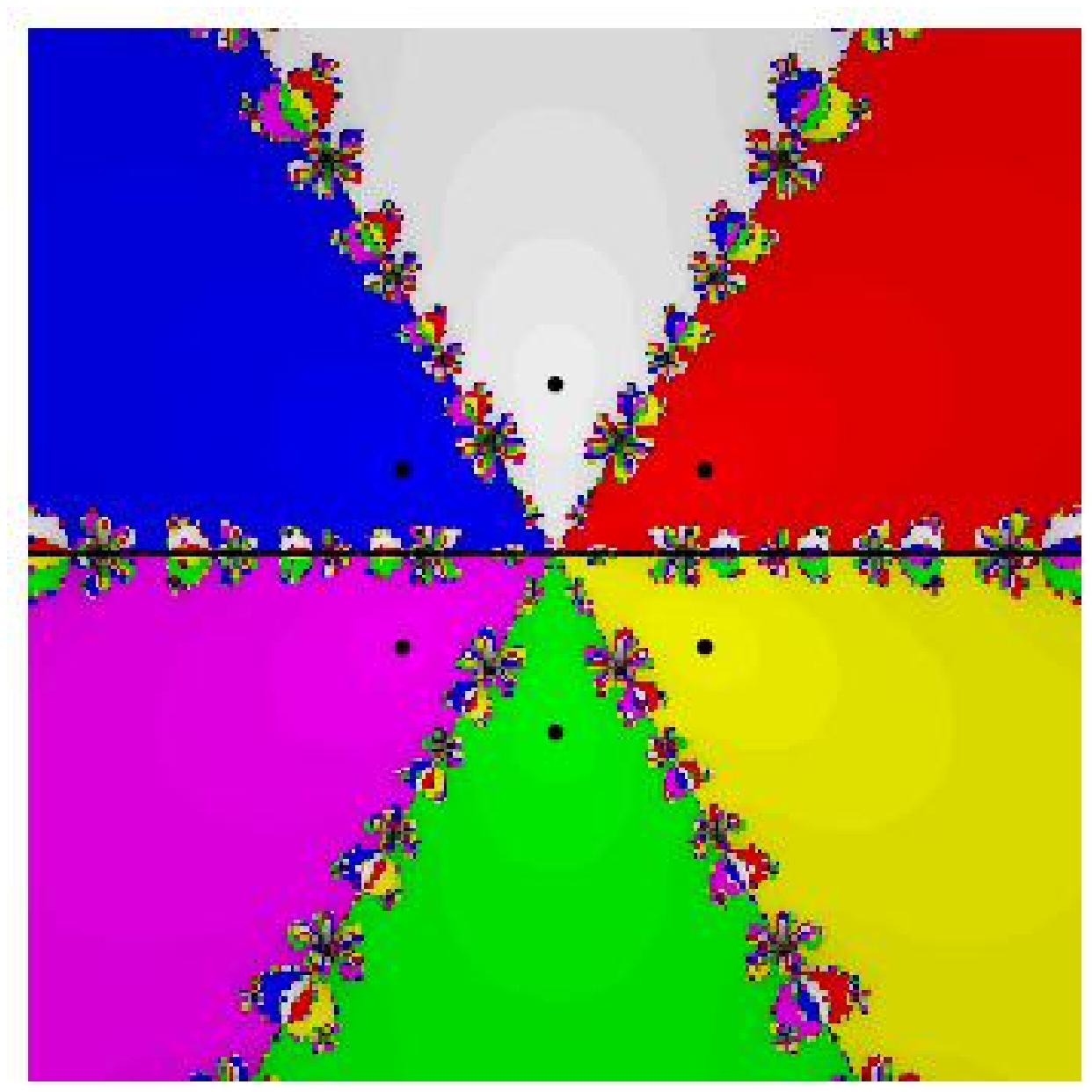}
\end{minipage}
\caption{Method (\ref{m2}) for test problems
$p_1(z)=\left(z+\frac{1}{z}\right)^5$, $p_2(z)=(z^3+1)^3$ and
$p_3(z)=\left(z+\frac{1}{z}\right)^2$ respectively}
\label{fig:figure4}
\end{figure}

\begin{figure}[ht!]
\begin{minipage}[b]{0.28\linewidth}
\centering
\includegraphics[width=\textwidth]{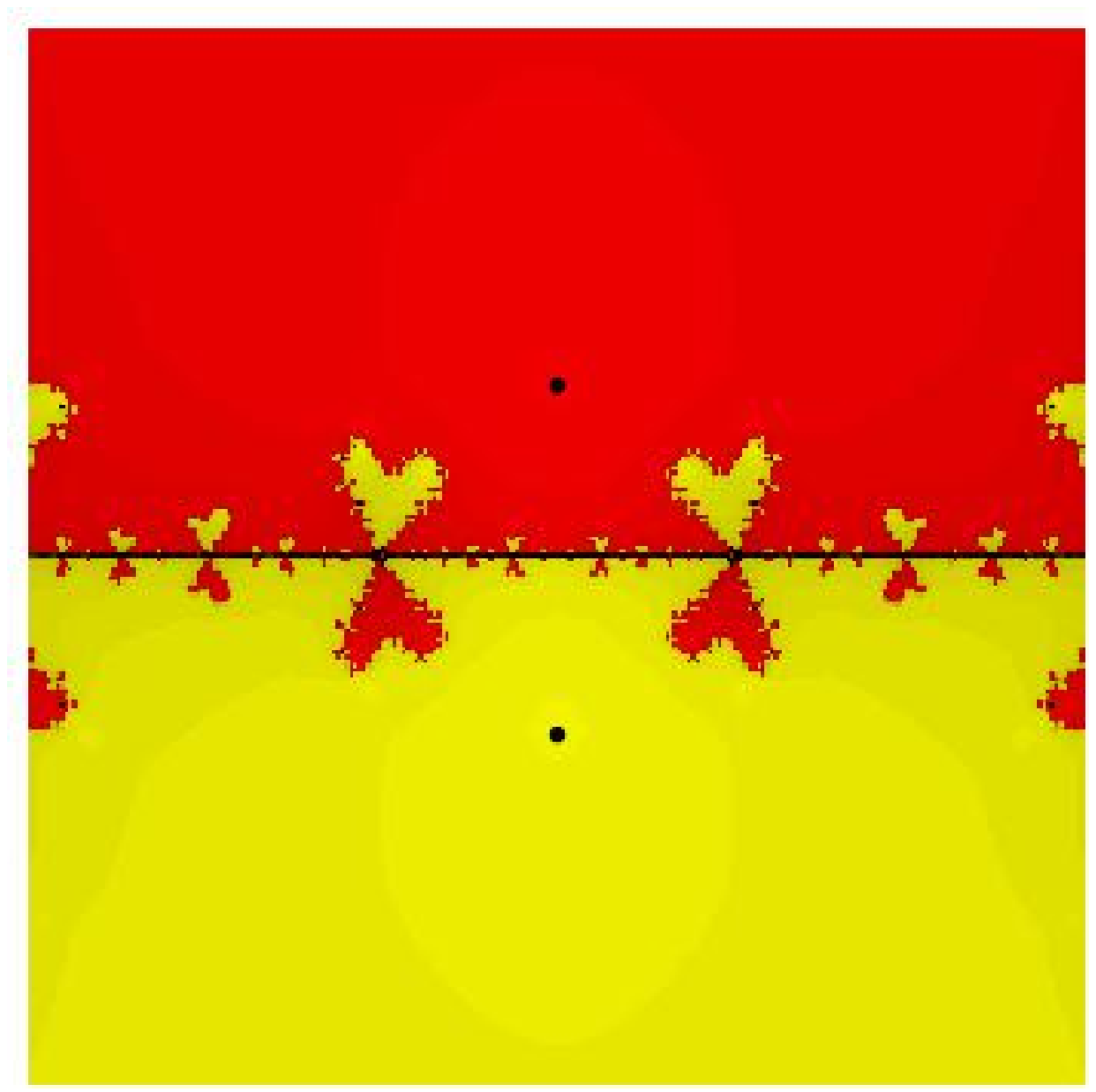}
\end{minipage}
\hspace{0.5cm}
\begin{minipage}[b]{0.28\linewidth}
\centering
\includegraphics[width=\textwidth]{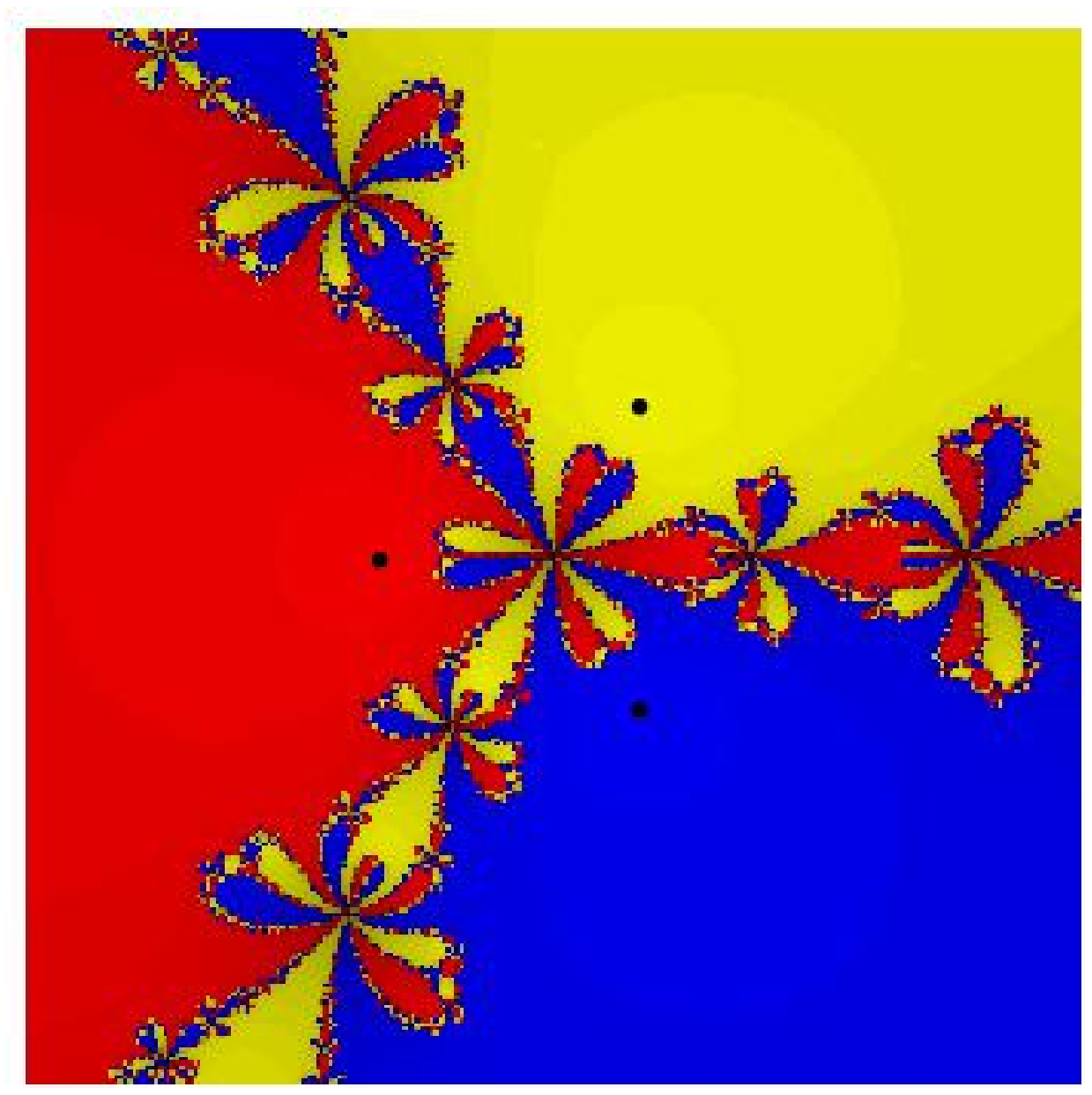}
\end{minipage}
\hspace{0.5cm}
\begin{minipage}[b]{0.28\linewidth}
\centering
\includegraphics[width=\textwidth]{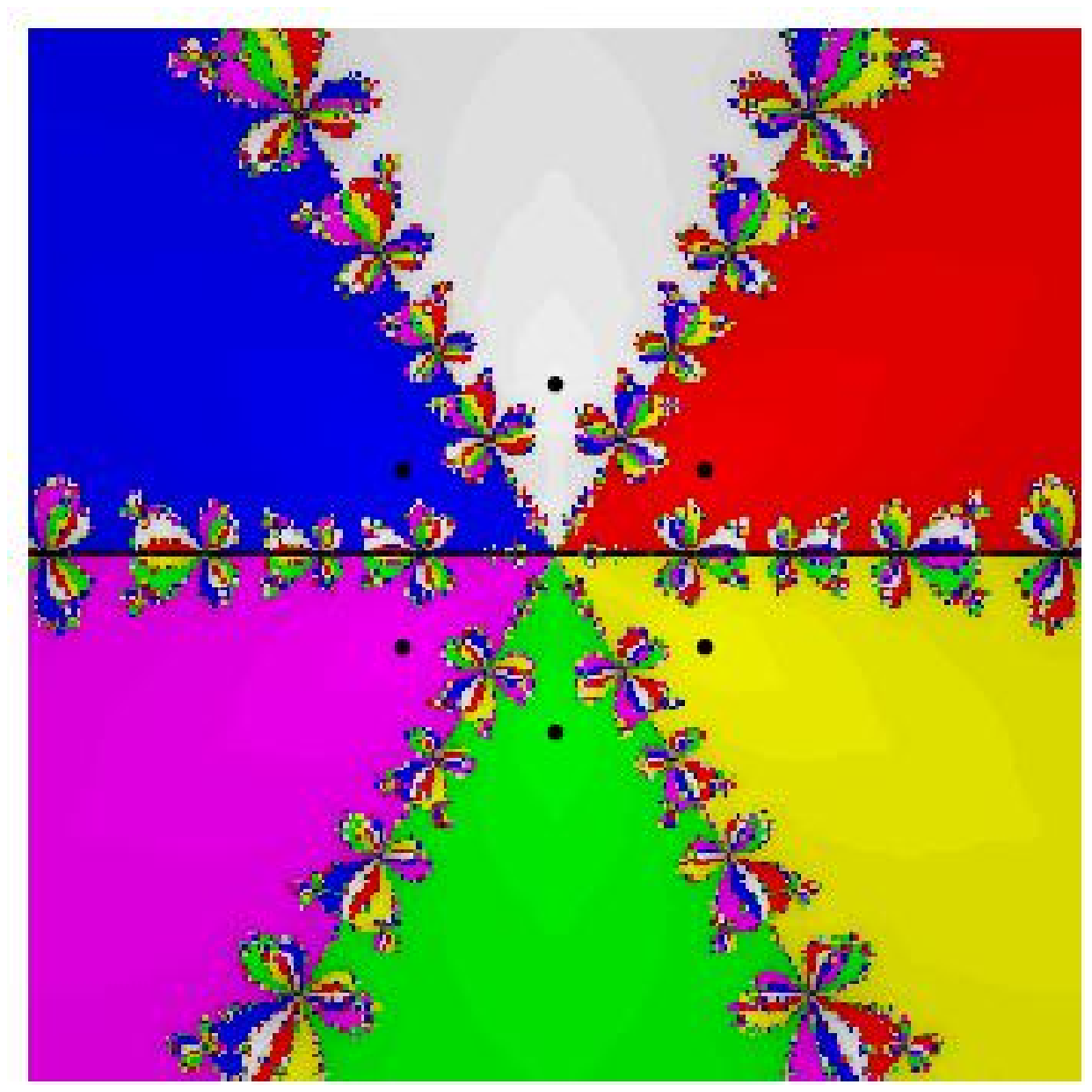}
\end{minipage}
\caption{Method (\ref{m3}) for test problems
$p_1(z)=\left(z+\frac{1}{z}\right)^5$, $p_2(z)=(z^3+1)^3$ and
$p_3(z)=\left(z+\frac{1}{z}\right)^2$ respectively}
\label{fig:figure5}
\end{figure}
\newpage
In Figures \ref{fig:figure1}-\ref{fig:figure5}, basins of
attractions of methods (\ref{c1}), (\ref{c4}), (\ref{m1}),
(\ref{m2}) and (\ref{m3}) are illustrated for three test problems
$p_1(z)$, $p_2(z)$ and $p_3(z)$ from left to the right
respectively. Therefore, the results presented in Table
\ref{table2} and Figures \ref{fig:figure1}-\ref{fig:figure5} show
that our proposed modified Potra-Pt{\'a}k method is competitive to
other existing methods.

\section{Conclusion}
\label{sec:conclusion} We have obtained a new method based on
Potra-Pt{\'a}k method for approximating multiple roots of
non-linear equations with same order of convergence and without
any additional evaluations of the function or its derivatives. It
has an efficiency index to $3^{\frac{1}{3}}\approx 1.44225$. The
theoretical results have been checked with some numerical
examples. On the other hand, numerical examples as well as the
basin of attraction show that our method works and can compete
with other methods in the same class.

\end{document}